\documentclass{amsart}
\usepackage{amsmath,amssymb,amscd}
\usepackage{eucal}

\newtheorem{theorem}{Theorem}
\newtheorem{proposition}[theorem]{Proposition}

\newtheorem{lemma}[theorem]{Lemma}

\newenvironment{definition}{
{\medskip\par\noindent\bf Definition.
}}{\vskip 2ex\par}
\newenvironment{remark}{
{\medskip\par\noindent\bf Remark.}}{\vskip 2ex\par}
\newcommand{\ga}{\alpha}

\newcommand{\gi}{\iota}

\newcommand{\call}{\mathcal{L}}
\newcommand{\fraks}{\mathfrak{s}}

\newcommand{\comp}{\mathrel{\scriptstyle\circ}}

\newcommand{\difftt}{\left. \dfrac{d{ }}{dt}\right|_{t=0} }

\newcommand{\Ad}{\operatorname{Ad}}

\newcommand{\Basic}{\operatorname{Basic}}

\newcommand{\C}{\mathbb{C}}

\newcommand{\cinf}{C^{\infty}}

\newcommand{\g}{\mathfrak{g}}

\newcommand{\half}{\frac{1}{2}}

\newcommand{\inv}{^{-1}}

\newcommand{\ix}{\iota_{X}}

\newcommand{\Lie}{\mathcal{L}}
\newcommand{\Liealg}{\operatorname{Lie}}
\newcommand{\liex}{\Lie_{X}}

\newcommand{\pt}{\operatorname{pt}}

\newcommand{\R}{\mathbb{R}}

\newcommand{\separate}{\medskip}
\newcommand{\seq}[2]{#1_{1}, \ldots, #1_{#2}}

\newcommand{\sone}{S^{1}}

\newcommand{\stl}{S^{2r +1}}

\newcommand{\tx}{\tilde{X}}

\begin{document}

\title{Equivariant Characteristic Classes in the Cartan Model}

\dedicatory{Dedicated to the memory of V. K. Patodi}

\author{Raoul Bott}

\address{Department of Mathematics, Harvard University, Cambridge, MA 02138}
\email{bott@math.harvard.edu}
\thanks{This work and the first author were supported in 
part by the National Science Foundation.  The second author was supported in part by the 
Institut Henri Poincar\'e, the Ecole Normale Sup\'erieure, and the 
Minist\`ere de l'Enseignement Sup\'erieur et de la Recherche, France.}

\author{Loring W. Tu}

\address{Department of Mathematics, Tufts University, Medford, MA 02155-7049}
\email{ltu@tufts.edu}



\date{January 31, 2001}

\begin{abstract}
This note shows the compatibility of the differential geometric 
 and the topological formulations of equivariant characteristic 
classes for a compact connected Lie group action.
\end{abstract}

\maketitle

Suppose $G$ and $S$ are two compact Lie groups, and $P$ and $M$ are manifolds.  
A principal $G$-bundle $\pi: P \to 
M$ is said to be $S$-{\it equivariant} if $S$ acts on the left on both $P$ and 
$M$ in such a way that
\begin{enumerate}
\item[a)] the projection map $\pi$ is $S$-equivariant:
\[
\pi(s.p)=s.\pi (p) \quad \mbox{for all } s \in S \ \mbox{and } p\in P;
\]
\item[b)] the left action of $S$ commutes with the right action of $G$:
\[
s.(p.g)=(s.p).g  \quad \mbox{for all\ } s \in S, p\in P, \ \mbox{and } g\in 
G.
\]
\end{enumerate}

An $S$-equivariant principal $G$-bundle $\pi: P\to M$ 
induces a principal $G$-bundle $\pi_S: P_S \to M_S$ over the homotopy 
quotient $M_S$.  
In the topological category, the \emph{equivariant characteristic classes} 
of $P \to M$ are defined to be the corresponding 
ordinary characteristic classes of $P_S \to M_S$.  
Thus, the equivariant characteristic classes are elements of the equivariant 
cohomology ring $H_S^*(M)$.

There is also a differential geometric definition of equivariant 
characteristic classes in terms of the curvature of a connection on $P$ 
(\cite{berline-vergne})(\cite{berline-vergne83}).  
However, there does not seem to be an explanation or justification in the 
literature bridging the two approaches.
The modest purpose of this note is to show the compatibility of the usual 
differential geometric formulation of equivariant characteristic classes 
with the topological one.  We have also tried to be as self-contained as 
possible, which partly explains the length of this article.

Let us first recall the situation for ordinary characteristic classes.  
Here the famous Chern-Weil construction represents
the ordinary characteristic classes of a 
principal $G$-bundle $\pi:P\to M$ by differential forms as follows.  
Fix a connection 
$\Theta$ on $P$.  Then its curvature $K$ is a $2$-form on $P$ with values 
in the Lie algebra $\g$ of $G$.  For each $\Ad G$-invariant polynomial $f$ 
on $\g$, one shows that $f(K)$ is a basic form on $P$ and so is the 
pullback $\pi^*\Lambda_f$ of a form $\Lambda_f$ on $M$.  Moreover, 
$\Lambda_f$ is closed and the cohomology class of $\Lambda_f$ in $H^*(M)$ is 
independent of the connection.  The characteristic class $c_f(P)$ 
of $P$ associated to $f$ is the cohomology class $[\Lambda_f] \in H^*(M)$.

For an $S$-equivariant principal $G$-bundle $P\to M$, the equivariant 
characteristic classes live in the equivariant cohomology ring $H_S^*(M)$.  
Since the equivariant cohomology of $M$ is the cohomology of its Cartan 
model, it is natural to ask if equivariant characteristic classes  can be 
constructed in the Cartan model out of the curvature of a connection on $P$.

This is indeed possible.  Assume that the Lie group $S$ is compact connected
with Lie 
algebra $\fraks$.  Then there exists an $S$-invariant connection $\Theta$ 
on the $S$-equivariant principal $G$-bundle $\pi: P\to M$.  Let $K$ be the 
curvature form of the $S$-invariant connection $\Theta$.  For each $X 
\in \fraks$, define $L_X: P\to \g$ by $L_X= - \iota_X\Theta$.  We will show 
that the equivariant characteristic class associated to an $\Ad (G)$-invariant polynomial $f$ 
on $\g$ is represented by $f(K+L_X)$ in the Cartan model.
In terms of a basis $X_1, \dots, X_{\ell}$ for $\fraks$ and dual basis 
$\seq{u}{\ell}$ for $\fraks^*$, this element is $f(K+\sum u_i L_{X_i})$.

The outline of the proof is as follows.  Since the total space $ES$ of the 
universal $S$-bundle is not a manifold, the first obstacle in a 
differential geometric treatment of equivariant characteristic classes 
is that neither $P_S$ nor $M_S$ are manifolds.  
Nonetheless, by approximating $ES$ by finite-dimensional manifolds 
$ES(r)$, $r=1, 2, \dots$, we can approximate the homotopy quotients $P_S$ 
and $M_S$ by the manifolds
\[
P^r := ES(r) \times_S P, \qquad M^r := ES(r) \times_S M.
\]
The natural map
\begin{align*}
\tau: P^r &\to M^r \\
[x,p] &\mapsto [x, \pi(p)], \qquad x \in ES(r), \ p\in P
\end{align*}
is again a principal $G$-bundle and is a finite-dimensional approximation 
to $\pi_S :P_S \to M_S$.

On the principal $G$-bundle $ES(r)\times P \to ES(r)\times M$ we are able 
to write down a connection which is basic with respect to the $S$-action.  
This connection then descends to a connection $\Xi_r$ on $P^r \to M^r$.  
By computing explicitly the curvature $K_r$ of the connection $\Xi_r$, it 
is easy to determine the element of the Weil model of $H_S^*(M)$ that 
restricts to the characteristic classes $[f(K_r)] \in H^*(M^r)$ for all 
$r$.  Using the Mathai-Quillen isomorphism between the Weil model and the 
Cartan model (\cite{mathai-quillen}), we can then write down the element of the Cartan model that 
represents the equivariant characteristic class $c_f(P)$.

\separate
\noindent
{\bf Notations.}  Since there are two groups acting on $P$, we need to be 
careful to distinguish between them.  A general element of $S$ will be 
denoted $s$, and a general element of $G$ will be denoted $g$.  Let 
${\ell}_s$ and $r_g$ denote 
left multiplication by $s$ and right multiplication by $g$, respectively.  The 
dimension of $S$ will be $\ell$ and the dimension of $G$ will be $n$.  A basis 
for the Lie algebra $\fraks$ of $S$ is $\seq{X}{\ell}$, and a basis for the 
Lie algebra $\g$ of $G$ is $\seq{e}{n}$.  We denote by $\tilde{Y} = 
\tilde{Y}_P$ the fundamental vector field on $P$ corresponding to $Y \in 
\fraks$ under the $S$-action, and by $Z^*=Z_P^*$ the fundamental vector 
field on $P$ corresponding to $Z\in \g$ under the $G$-action.

\section{Equivariant Cohomology}

We begin with a rapid review of the basic constructions of equivariant 
cohomology.

Suppose $M$ is a CW-complex on which a Lie group $S$ acts on the left.  We 
call such a space an \emph{$S$-space}.  Let $\pi: ES \to BS$ be the 
universal $S$-bundle over the classifying space $BS$.  Then the total 
space $ES$ is a contractible space on which $S$ acts freely on the right.  
Thus, $ES \times M$ is a space with the same homotopy type as $M$, and $S$ 
acts freely on $ES\times M$ by the diagonal action:
\[
(e,x).s = (es, s\inv x), \qquad \mbox{for } e\in E, x\in M, s\in S.
\]
The quotient space $(ES\times M)/S$, denoted $ES \times_S M$ or $M_S$, is 
called the \emph{homotopy quotient} of $M$ by $S$.

\begin{definition}
The \emph{equivariant cohomology} $H_S^*(M)$ of $M$ under the action of 
$S$ is the singular cohomology of the homotopy quotient $M_S$.
\end{definition}

One may compute equivariant cohomology with any coefficient ring.  Since 
we are interested in manifolds and differential forms, we will take 
the coefficient ring to be $\R$. 

The homotopy quotient of a point is the classifying space of $S$, since
\[
(\pt)_S = (ES\times \pt )/S = (ES)/S = BS.
\]
Hence, the equivariant cohomology of a point is the ordinary cohomology of 
the classifying space $BS$.

A map $f: M \to N$ of $S$-spaces is \emph{$S$-equivariant} or \emph{an 
$S$-map} if $f(s.x)=s.f(x)$ for all $x$ in $M$.  An $S$-equivariant map 
$f:M\to N$ induces a map $f_S: M_S \to N_S$ of homotopy quotients and 
therefore a homomorphism in equivariant cohomology:
\[
f_S^*: H_S^*(N) \to H_S^*(M).
\]
Since every $S$-space $M$ has an $S$-equivariant map $f: M \to \pt$, there 
is a ring homomorphism
\[
f_S^*: H_S^*(\pt) \to H_S^*(M).
\]
Thus, the equivariant cohomology ring $H_S^*(M)$ is a module over $H^*(BS)$.

In summary, for a fixed Lie group $S$, equivariant cohomology is a 
contravariant functor from the category of $S$-spaces and $S$-maps to the 
category whose objects are rings that are simultaneously 
$H^*(BS)$-modules and whose morphisms are ring homomorphisms that are 
simultaneously $H^*(BS)$-module homomorphisms.

\section{The Weil Model}

As before, let $S$ be a Lie group acting on a space $M$.  When $M$ is a 
manifold and the action is $\cinf$, we might ask whether it is possible to 
compute equivariant cohomology using differential forms, just as the de Rham 
complex of $M$ computes its singular cohomology with real coefficients.  
The Weil model and the Cartan model answer this question affirmatively.

Let $X_1, \dots, X_{\ell}$ be a basis for the Lie algebra $\fraks$ of the Lie 
group $S$.  The Lie algebra structure of $\fraks$ is completely specified 
by its structure constants $c_{ij}^k$, defined by
\[
[X_i,X_j]=\sum_{k=1}^{\ell} c_{ij}^k X_k.
\]

Suppose that a Lie group $S$ acts on a manifold $P$.  Then each $Y$ in the 
Lie algebra $\fraks$ defines a \emph{fundamental vector field} $\tilde{Y}=
\tilde{Y}_P$ on 
$P$ as follows:  for $p\in P$,
\[
\tilde{Y}_p = 
\begin{cases}
\difftt e^{-tY}.p &\mbox{if $S$ acts on $P$ on the left;}\\
\difftt p.e^{tY} &\mbox{if $S$ acts on $P$ on the right.}
\end{cases}
\]

If $S$ acts on the left on $P$ and $i_p: S \to P$ is the map $i_p(s)=s.p$, 
then $\tilde{Y}_p= (i_p)_*(-Y)$; similarly, if $S$ acts on the right on 
$P$ and $i_p(s)=p.s$, then $\tilde{Y}_p= (i_p)_*(Y)$.  
Our choice of signs ensures that the 
construction of fundamental vector fields is a Lie algebra homomorphism,
whether the group $S$ acts on the left or on the right: 
for $X,Y \in \fraks$,
\begin{equation} \label{e:liehomo}
[X,Y]\,\tilde{} =[\tilde{X}, \tilde{Y}].
\end{equation}

Recall that a \emph{connection} on a principal $S$-bundle $\pi: 
P \to M$ is an $\fraks$-valued $1$-form $\omega$ on $P$ such that
\begin{enumerate}
\item[a)] if $\tilde{Y}_P$ is the fundamental vector field on $P$ associated
to $Y \in \fraks$, then $\omega(\tilde{Y}_P)=Y$, and
\item[b)] if $r_s$ is right translation by $s \in S$, then $r_s^*\omega= (\Ad 
s\inv)\omega$.
\end{enumerate}

In terms of a basis $\seq{X}{\ell}$ for the Lie algebra $\fraks$,
\[
\omega= \sum \omega_i X_i,
\]
where $\seq{\omega}{\ell}$ are real-valued $1$-forms on $P$.  In differential 
geometry one shows that the connection forms $\omega_i$ and the curvature 
forms $\Omega_i$ satisfy the equations
\begin{enumerate}
\item[1)] $d\omega_i + \frac{1}{2} \sum_{j,k}\, c_{jk}^i \omega_j\wedge 
\omega_k = \Omega_i$,
\item[2)] $d\Omega_i = \sum_{j,k}\, c_{jk}^i \Omega_j \wedge \omega_k$.
\end{enumerate}

For differential forms $\omega= \sum \omega_i X_i$ and $\tau = \sum \tau_j 
X_j$ with values in the Lie algebra $\fraks$, we introduce the notation
\[
[\omega, \tau]= \sum \omega_i \wedge \tau_j [X_i, X_j].
\]
Then we may rewrite 1) and 2) as
\begin{enumerate}
\item[$1'$)] $d\omega + \frac{1}{2} [\omega, \omega] = \Omega$,
\item[$2'$)] $d\Omega = [\Omega, \omega]$.
\end{enumerate}

The \emph{Weil algebra} $W(\fraks)$ of the Lie algebra $\fraks$ is defined 
to be the tensor product of the symmetric algebra and the exterior algebra
of $\fraks^*$:
\[
W(\fraks)= S^*(\fraks^*) \otimes \Lambda^*(\fraks^*).
\]
Fix a basis $\seq{X}{\ell}$ for $\fraks$.  To distinguish the elements of 
$S^*(\fraks^*)$ from those in $\Lambda^*(\fraks^*)$, we will write the 
dual basis for $\fraks^*$ in $S^*(\fraks^*)$ as $\seq{u}{\ell}$, and the dual 
basis for $\fraks^*$ in $\Lambda^*(\fraks^*)$ as $\seq{\theta}{\ell}$.  We give 
the Weil algebra a grading by
\[
\deg u_i = 2, \qquad \deg \theta_i = 1.
\]

Define the \emph{Weil differential} on $W(\fraks)$ to be the 
antiderivation $d$ such that
\[
d\theta_i + \frac{1}{2} \sum_{j,k} c_{jk}^i \theta_j\theta_k= u_i, \qquad 
du_i = \sum_{j,k} c_{jk}^i u_j\theta_k.
\]
For each $X\in \fraks$, define an antiderivation $\gi_X$ on $W(\fraks)$ by
\[
\gi_X\theta_i = \theta_i (X), \qquad \gi_X u_i = 0.
\]
Finally, define the Lie derivative $\liex$ on $W(\fraks)$ by the formula
\[
\liex = d\gi_X + \gi_X d.
\]

Because the de Rham complex $\Omega^*(M)$ of the $S$-manifold $M$ also has 
these operators $d$, $\gi_X$ and $\liex$, the tensor product 
$W(\fraks)\otimes \Omega^*(M)$ inherits the same operators in the usual 
way.  An element $\ga$ of $W(\fraks)\otimes \Omega^*(M)$ is said to be 
\emph{basic} if 
\[ 
\gi_X\ga=0 \quad \mbox{ and }\quad 
\liex \ga = 0.
\]
The first condition says $\ga$ is \emph{horizontal}, and the second says 
$\ga$ is \emph{invariant}.
It is easily shown that if $\ga$ is basic, then so is $d\ga$.  Hence, the 
space of basic forms, $\Basic (W(\fraks)\otimes \Omega^*(M))$, is a 
differential complex under $d$.

\separate
\begin{theorem}
Suppose a Lie group $S$ acts on the manifold $M$.  Then
\[
H^*\{\Basic(W(\fraks)\otimes \Omega^*(M))\}\simeq H_S^*(M).
\]
\end{theorem}

We call the differential complex $\Basic(W(\fraks)\otimes \Omega^*(M))$ the 
\emph{Weil model} for the equivariant cohomology of $M$ under the action of 
$S$.

\section{The Cartan Model}

An element of the Weil model is a finite sum
\[
\ga= a+\sum \theta_i a_i + \sum_{i<j} \theta_i\theta_j a_{ij} + \dots
+ \theta_1 \dots \theta_n a_{1\dots n},
\]
where $a, a_i, a_{ij}, \dots \in S^*(\fraks^*)\otimes \Omega^*(M)$.
We now write $\gi_i$ for $\gi_{X_i}$.  The horizontality condition implies 
that
\[
a_i = - \gi_i a, \quad a_{ij}= \gi_i a_j= -\gi_i\gi_j a,
\]
and in general,
\begin{equation} \label{e:coef}
a_{i_1 \dots i_k} = (-1)^k \gi_{i_1} a_{i_2 \dots i_k}.
\end{equation}
By induction, all the coefficients $a_I$ of a horizontal form $\ga$ are 
determined by $a \in S^*(\fraks^*)\otimes \Omega^*(M)$.

The Lie group $S$ acts on the symmetric algebra $S^*(\fraks^*)$ by the 
coadjoint representation and on $\Omega^*(M)$ by pulling back forms.  
Hence, $S$ also acts on the tensor product $S^*(\fraks^*)\otimes 
\Omega^*(M)$.  We denote the invariant subspace by
$(S^*(\fraks^*)\otimes \Omega^*(M))^S$.

It follows from (\ref{e:coef}) that every element of $\Basic(W(\fraks^*)\otimes \Omega^*(M))$ 
can be written in the form
\begin{align} \label{e:correspondence}
\ga &= a - \sum \theta_i \gi_i a + \sum_{i<j} (\theta_i\gi_i)(\theta_j 
\gi_j)a - \dots \\
&=(1-\theta_1 \gi_1) \dots (1-\theta_n\gi_n)a
\end{align}
for a unique $a\in (S^*(\fraks^*)\otimes \Omega^*(M))^S$.

\separate
\begin{theorem}
There is an isomorphism, called the Mathai-Quillen isomorphism,
\begin{align*}
\Basic(W(\fraks^*)\otimes \Omega^*(M)) &\leftrightarrow 
(S^*(\fraks^*)\otimes \Omega^*(M))^S, \\
\ga &\mapsto a \qquad \mbox{\rm (as in Eq. (\ref{e:correspondence}))} \\
(1-\theta_1 \gi_1) \dots (1-\theta_n\gi_n)a &\leftarrow a.
\end{align*}
Under this isomorphism, the Weil differential $d$ corresponds to the 
\emph{Cartan differential}
\[
d_C a = da - \sum u_i\gi_ia
\]
on $(S^*(\fraks^*)\otimes \Omega^*(M))^S$.
\end{theorem}

\separate

Thus, the differential complex $\{(S^*(\fraks^*)\otimes \Omega^*(M))^S, 
d_C\}$ also computes the equivariant cohomology $H_S^*(M)$.  It is called 
the \emph{Cartan model} of $H_S^*(M)$.

For a circle action on $M$, we fix an element $X_1=X \in \Liealg (S^1)$ 
and choose $\theta$ to be its dual $1$-form.  Since $S^1$ is abelian, its 
coadjoint representation on $S^*(\fraks^*)$ is trivial.  Hence, the Cartan 
model is $\Omega^*(M)^{S^1}[u]$.  The element of the Weil model
corresponding to the element $a$ in the Cartan model is
\[
\ga = a - \theta \gi_X a, \qquad a \in  \Omega^*(M)^{S^1}[u].
\]

\section{Fundamental Vector Fields of the $S$-Action}

Let $\pi: P\to M$ be an $S$-equivariant principal $G$-bundle.  By 
hypothesis, $S$ acts on the left on $P$, and $G$ acts on the right on $P$, 
and the two actions commute.

\separate  

\begin{proposition} \label{p:funvecfield}
For $X \in \fraks$, let 
$\tilde{X}=\tilde{X}_P$ be the fundamental vector field corresponding to 
$X$ under the $S$-action.  Then $\tilde{X}$ is invariant under $G$ and 
transforms by the adjoint representation under $S$; more precisely,
\begin{enumerate}
\item[i)] for $g\in G$, $r_{g*}(\tilde{X})=\tilde{X}$;
\item[ii)] for $s \in S$, ${\ell}_{s*}(\tilde{X})= ((\Ad s)X)\tilde{}$.
\end{enumerate}
\end{proposition}

\begin{proof}
i)  Let $i_p: S \to P$ be the map $i_p(s)=s.p$.  Then $r_g \comp i_p = 
i_{pg}$.  Hence,
\[
r_{g*}\tx_p = r_{g*}i_{p*}(-X)=i_{pg*}(-X) = \tx_{pg}.
\]
ii) For $s,t \in S$,
\[
{\ell}_s \comp i_p(t) = stp=sts\inv s p = i_{sp}\comp c_s (t).
\]
where $c_s (t)= sts\inv$ is conjugation by $s$.  Then
\begin{align*}
({\ell}_s)_* \tilde{X}_p &={\ell}_{s*} i_{p*} (-X) = i_{sp*}c_{s*}(-X)
=i_{sp*}(-(\Ad s) X) \\
&=((\Ad s)X)\tilde{}_{sp}. \phantom{xxxxxxxxxxxxxxxxxxxxxxxxxx}
\end{align*}
\end{proof}

If $\omega$ is a differential form on $P$ and $X$ is an element of the Lie 
algebra $\fraks$, we define $\iota_X \omega$ to be the contraction 
$\iota_{\tilde{X}} \omega$, and $\call_X \omega$ to be the Lie derivative 
$\call_{\tilde{X}}\omega$.

\section{The $\g$-Valued Function $L_X$}

Recall that a \emph{connection} on the principal $G$-bundle $P\to M$ is 
a $\g$-valued $1$-form $\Theta$ on $P$ satisfying the two conditions:
\begin{enumerate}
\item[i)] $\Theta(Z_P^*)=Z$ if $Z\in\g$ and $Z_P^*$ 
is the fundamental vector 
field on $P$ associated to $Z$ by the $G$-action,
\item[ii)] $r_g^*\Theta = (\Ad g\inv)\Theta$ for all $g\in G$.
\end{enumerate}

Fix a connection $\Theta$ on the $S$-equivariant principal $G$-bundle
$\pi:P\to M$.  
For $X \in \fraks$, define the function $L_X : P \to \g$ by
\[
L_X (p) := -(\iota_X \Theta)(p):= -\Theta_p (\tilde{X}_p).
\]

When the principal $G$-bundle $P \to M$ is $S$-equivariant, we say that 
the connection $\Theta$ is \emph{$S$-invariant} if the Lie derivative 
$\Lie_{X}\Theta = 0$ for all $X \in \fraks$.
 This is equivalent to 
${\ell}_s^*\Theta=\Theta$ for all $s\in S$.
The compactness of $S$ allows us to average any connection $\Theta$ over 
$S$ to obtain an $S$-invariant connection
\[
\int_S {\ell}_s^* \Theta \,ds.
\]
Since characteristic classes are independent of the connection, we may 
assume from the outset, if we like, 
that the connection $\Theta$ is $S$-invariant. 

\separate 

\begin{proposition} \label{p:translatelx}
\begin{enumerate}
\item[i)] If $\Theta$ is any connection on $P$, then
for $g\in G$, $r_g^* L_X = (\Ad g\inv)L_X$.
\item[ii)] If the connection $\Theta$ is $S$-invariant, then for $s \in 
\fraks$, $\ell_s^* L_X = L_{(\Ad s\inv)X}$.
\end{enumerate}
\end{proposition}

\begin{proof}
i) By Lemma \ref{p:funvecfield}(i), the vector field $\tx$ is 
right-invariant under $G$. Thus, for $p\in P$ and $g\in G$,
\begin{align*}
L_X (pg) &=-\Theta_{pg}(\tx_{pg}) = -\Theta_{pg} (r_{g*}\tx_p) 
=-(r_g^*\Theta_{pg} )(\tx_p) \\
&= -(\Ad g\inv)\Theta_{p} (\tx_p) = (\Ad g\inv) L_X(p).
\end{align*} 
\separate

\noindent
ii) 
\begin{alignat*}{2}
({\ell}_s^*L_X)(p)= L_X (sp) &= - \Theta_{sp} (\tilde{X}_{sp})&& \\
&=-\Theta_{sp} ({\ell}_{s*}((\Ad s\inv)X)\tilde{}_p) &&\quad\text{(by 
Prop.~\ref{p:funvecfield}(ii))} \\
&=-{\ell}_s^*\Theta_{sp}(((\Ad s\inv)X)\tilde{}_p) &&\\
&=-\Theta_p (((\Ad s\inv)X)\tilde{}_p) &&\quad\text{(by the $S$-invariance 
of $\Theta$)} \\
&= L_{(\Ad s\inv)X}. &&\phantom{by the S-invariance of Theta} 
\end{alignat*} 
\end{proof}

Although the connection $\Theta$ is $S$-invariant, Prop.~\ref{p:translatelx}(ii) shows 
that the function $L_X$ is not $S$-invariant.  Indeed, on the Lie algebra 
level, we have the following identity.

\separate

\begin{proposition} \label{p:lielx}
If $\Theta$ is an $S$-invariant connection on $P$ and $X,Y \in \fraks$, 
then $\call_Y L_X = L_{[Y,X]}$.
\end{proposition}

\begin{proof}
\begin{alignat*}{2}
\call_Y L_X &= -\call_Y \Theta(\tilde{X}) && \\
&=-(\call_Y\Theta)(\tilde{X}) - \Theta(\call_Y \tilde{X}) & 
&\quad\text{(derivation property of the Lie derivative)} \\
&=-\Theta([\tilde{Y},\tilde{X}]) &&\quad\text{(by the $S$-invariance of 
$\Theta$)} \\
&=-\Theta([Y,X]\,\tilde{}) &&\quad\text{(by Eq.~(\ref{e:liehomo}))} \\
&= L_{[Y,X]}. && \phantom{xxxxxxxxxxxxxxxxxxxxxxxxxxx} 
\end{alignat*}
\end{proof}

If $X_1, \dots, X_{\ell}$ is a basis for the Lie algebra $\fraks$,
we write $\iota_i$ and $L_i$ for $\iota_{X_i}$ and $L_{X_i}$, respectively.

\section{The Curvature of a Connection}

The curvature of a connection 
$\Theta$ on a principal $G$-bundle $P\to M$ is given by the well-known 
formula:
\[
K=d\Theta +\frac{1}{2}[\Theta,\Theta].
\]
In case $G$ is a matrix group, this is equivalent to
\[
K=d\Theta + \Theta \wedge \Theta.
\]
By the equivariance property of a connection, it is easily checked that 
under right translation $r_g$ by an element $g$ in $G$,
\[
r_g^*K=(\Ad g\inv)K.
\]
In the terminology of Kobayashi and Nomizu (\cite{kobayashi-nomizu}), 
this shows that $K$ is a 
``pseudotensorial $2$-form of type $\Ad G$", and hence
\[
K\in \Omega^2(M,\Ad P).
\]

In general, we can identify $\Omega^k (M, \Ad P)$ with the space of 
horizontal $\g$-valued 
$k$-forms $\ga$ on $P$ satisfying
\[
r_g^* \ga = (\Ad g\inv) \ga.
\]
Thus, for any connection $\Theta$ on an $S$-equivariant principal 
bundle $\pi: P \to M$ and for any $X \in \fraks$,
\[
L_X \in \Omega^0(M, \Ad P), \qquad \Theta \in \Omega^1(M,\Ad P),
\qquad K\in \Omega^2(M,\Ad P).
\]

A connection on a principal bundle induces a covariant derivative on all 
associated bundles.  In particular, the connection $\Theta$ on the 
principal bundle $P\to M$ induces a covariant derivative $D$ on the adjoint 
bundle:
\[
D: \Omega^k (M, \Ad P) \to \Omega^{k+1} (M, \Ad P), \quad \mbox{for all } 
k= 0,1,2,\ldots .
\]
It is known that the covariant derivative $D$  on $\Omega^k (M, \Ad P)$
is given by the formula
\[
D \ga= d \ga + [\Theta, \ga].
\]

\separate

\begin{proposition} \label{p:ixk}
Let $P\to M$ be an $S$-equivariant principal $G$-bundle, $\Theta$ an 
$S$-invariant connection on $P$, and $K$ the curvature form of 
$\Theta$.  For $X,Y \in \fraks$,
\begin{enumerate}
\item[i)]
$\ix K = D L_X = dL_X +  [\Theta, L_X] $.
\item[ii)] $\iota_Y\iota_X K= [L_X,L_Y]-L_{[X,Y]}$.
\end{enumerate}
\end{proposition}

\begin{proof}
i) Applying $\ix$ to both sides of
\[
K=d\Theta +\frac{1}{2}[\Theta, \Theta],
\]
we get
\begin{align*}
\ix K &= \ix d\Theta + \frac{1}{2} [\ix\Theta, \Theta] - \frac{1}{2} 
[\Theta,\ix \Theta] &&\\
&=\liex \Theta -d\ix \Theta -
\frac{1}{2} [L_X, \Theta] +\frac{1}{2} [ \Theta, L_X] &&
 (\mbox{because } \liex = \ix d + d \ix) \\
&= dL_X + [\Theta, L_X] && (\liex \Theta = 0 \mbox{\ since $\Theta$ is 
$S$-invariant}).
\end{align*}
This last expression is precisely the covariant derivative $DL_X$ of $L_X$ 
induced by the connection $\Theta$, where we view $D$ as a map from 
$\Omega^0 (M, \Ad P)$ to $\Omega^1 (M, \Ad P)$. 

\separate
\noindent
ii) \begin{align*}
\iota_Y\iota_X K &= \iota_Y (dL_X +
 [\Theta, L_X]) && \text{(by Part (i))}\\
&=(\call_Y - d\iota_Y)L_X + [\iota_Y\Theta, L_X]  &&\\
&=\call_Y L_X -  [L_Y, L_X] && \\
&=L_{[Y,X]} +[L_X,L_Y] &&\text{(Prop.\ \ref{p:lielx})} \\
&=[L_X,L_Y]-L_{[X,Y]}.  &&\phantom{xxxxxxxxxxxxxxxxxxxxxxx} 
\end{align*}
\end{proof}

\section{The Circle Case}

At this point, we specialize to a circle action, not because of any 
logical necessity, but because when $S$ is the circle $S^1$, 
the computation is simpler and the ideas of the proof are more 
transparent.  In the circle case, the approximating bundle $ES(r) \to 
BS(r)$ is the Hopf fibration $S^{2r+1}\to \C P^r$.

We view the sphere $S^{2r+1}$ as the unit sphere in $\C^{r +1}$.  
Then $S^1$ acts on $S^{2r+1}$ on the right by scalar multiplication:
\[
z.t=(z_0, \ldots, z_{r}).t=(tz_0, \ldots, tz_{r}), \quad t \in \sone, 
z\in \stl.
\]
Clearly, the stabilizer at each point $z \in \stl$ is the identity, so this 
is a free action.
Fix a nonzero element $X$ in $\fraks = \Liealg (\sone)$ and let  
$\tilde{X} := \tilde{X}_{\stl}$ be the fundamental vector field on 
$\stl$ associated to $X$.
Let $\theta$ be a connection for the $\sone$-bundle $\stl \to \C P^r$.  
Since $\sone$ is abelian, 
\[
r_s^* \theta = (\Ad s\inv) \theta = \theta.
\]
So the connection $\theta$ is $S$-invariant.  In particular,
\begin{equation} \label{e:lxtheta}
\call_X \theta = 0.
\end{equation}
Moreover, since fixing $X$ amounts to choosing an isomorphism $\Liealg 
(\sone) \simeq \R$ under which $X$ corresponds to $1$,
\[
\iota_X \theta = \theta (\tilde{X}) \equiv 1.
\]

\subsection{Induced Connection and Curvature on $P^{r} \rightarrow M^{r}$}

Let $P\to M$ be an $\sone$-equivariant principal $G$-bundle with an 
$\sone$-invariant connection $\Theta$.  In the circle case, we fix an 
element $X \in \Liealg (\sone)$ and abbreviate $L_X := -\iota_X \Theta$ as $L$.

We approximate the homotopy quotients $P_{\sone}$ and $M_{\sone}$ by
\[
P^r := \stl \times_{\sone} P, \qquad M^r := \stl \times_{\sone} M.
\]
Then the principal $G$-bundle 
$P^{r} \to M^{r}$ is a finite-dimensional approximation to 
$P_{\sone} \to M_{\sone}$.  
Instead of finding directly a 
connection on  $P^{r} \to M^{r}$, we will construct a connection on 
$\stl\times P \to \stl \times M$ that is basic with respect to the 
$\sone$-action.  It will then descend to a connection on 
$P^{r} \to M^{r}$.

\separate
\begin{proposition}
The $\g$-valued $1$-form
\[
\Xi =  1\otimes \Theta + \theta \otimes L
\]
on $\stl\times P$ is a connection form for the principal $G$-bundle $\stl 
\times P \to \stl \times M$ and is basic with respect to the circle action 
on $\stl \times P$.
\end{proposition}

\begin{proof}
    Since $L$ and $\Theta$ are $\g$-valued and $\theta$ and $1$ are 
    scalar-valued, $\Xi$ is a $\g$-valued $1$-form on $\stl \times P$.
\begin{enumerate}
    \item[i)] If $Z \in \g$ and $Z_{P}^*$ is the fundamental vector 
    field on $P$ associated to $Z$ under the $G$-action, 
then $(0, Z_{P}^*)$ is the 
    fundamental vector field on $\stl \times P$ associated to $Z$.  
    Hence,
    \[
    \Xi (0, Z_{P}^*) = ( 1 \otimes \Theta + \theta \otimes L)(0, Z_{P}^*)
    = \Theta (Z_{P}^*) = Z.
    \]
    \item[ii)] Right-equivariance:  For $g\in G$,
    \begin{align*}
	r_{g}^{*}\Xi &= 1\otimes r_g^* \Theta + ( r_g^*\theta) \otimes r_g^* L\\
&= 1 \otimes (\Ad g\inv) \Theta + \theta\otimes (\Ad g\inv)L \\
&= (\Ad g\inv) ( 1 \otimes \Theta + \theta\otimes L) =(\Ad g\inv) \Xi.
\end{align*}
Hence, $\Xi$ is a connection on $\stl\times P \to \stl \times M$.
\end{enumerate}
To prove that $\Xi$ is basic with respect to the $\sone$-action, we check the horizontality condition $\ix 
\Xi=0$ and the invariance condition $\Lie_X\Xi = 0$.
\begin{enumerate}
\item[iii)] Horizontality:
\begin{align*}
\ix \Xi &= \ix ( 1 \otimes \Theta + \theta\otimes L) \\
&=1 \otimes \ix \Theta + (\ix \theta)\otimes L - \theta \otimes \ix L\\
&= -1\otimes L + 1\otimes L = 0.
\end{align*}
\item[iv)] Invariance:
\begin{align*}
\Lie _X \Xi &= \Lie_X ( 1\otimes \Theta +\theta\otimes L) \\
&= 1 \otimes \Lie_X \Theta + 
(\Lie_X \theta)\otimes L + \theta \otimes \Lie_X L.
\end{align*}
\end{enumerate}
In this sum, $\Lie_X \theta = 0$ and
$\Lie_X \Theta = 0$ because both connections $\theta$ and $\Theta$ are 
$\sone$-invariant (Eq.~\eqref{e:lxtheta}); moreover, by Prop.\ 
\ref{p:lielx},
\[
\Lie_X L = \Lie_X L_X = L_{[X,X]} = 0.
\]
Hence, $\Lie_X \Xi = 0$.
\end{proof}

So $\Xi=1\otimes\Theta+\theta\otimes L$ descends to $P^{r}=\stl 
\times_{\sone} P$, i.e., if $\ga: \stl \times P \to P^{r}$ is the 
projection, then 
$\Xi = \ga^* \Xi '$ for some $\g$-valued $1$-form $\Xi'$ on 
$P^{r}$.  From the following lemma, it follows that $\Xi'$ is a 
connection for $\tau:P^{r} \to M^{r}$.

\separate
\begin{lemma} \label{p:pullback}
Let $G$ be a Lie group with Lie algebra $\g$ and $\pi:P\to M$, $\pi':P'\to 
M'$ two principal $G$-bundles.  Suppose there is a $G$-equivariant 
submersion $\ga:P\to P'$ and $\Theta'$ is a $\g$-valued 1-form on $P'$.  
Then $\ga^*\Theta'$ is a connection for $\pi:P\to M$ if and only if 
$\Theta'$ is a connection for $\pi':P' \to M'$.
\end{lemma}

\separate

The proof of this lemma is based on a series of simple remarks.

\begin{remark}
\begin{enumerate}
\item[i)] For $p \in P$, let $i_p:G \to P$ be the map $i_p(g)=pg$.  
Similarly, for $p' \in P'$, let $i_{p'}: G \to P'$ be the map 
$i_{p'}(g)=p'g$.  Then 
\[
\ga \comp i_p = i_{\ga(p)}.
\]

\item[ii)] For $g \in G$, let $r_g$ be right translation by $g$ on a 
principal $G$-bundle.  Since $\ga:P \to P'$ is $G$-equivariant,
\[
r_g \comp \ga = \ga \comp r_g.
\]

\item[iii)] For $Y \in \g$, if $Y_P$ and $Y_{P'}$ are the fundamental 
vector fields associated to $Y$ on $P$ and $P'$ respectively, then 
$Y_{P'}= \ga_*(Y_P)$.
\end{enumerate}

\medskip
\noindent
{\it Proof of iii).}
 For $p\in P$,
\[
\ga_*(Y_P)_p = \ga_*{i_p}_*(Y)= (\ga\comp i_p)_* Y= {i_{\ga(p)}}_*Y=
(Y_{P'})_{\ga(p)}. 
\]
\end{remark}  
\smallskip

\noindent
\emph{Proof of Lemma \ref{p:pullback}.}
$(\Leftarrow)$ Suppose $\Theta'$ is a connection for $\pi':P'\to M'$, and 
$Y\in \g$.  Then
\begin{enumerate}
\item[a)] $(\ga^*\Theta')(Y_P)= \Theta' (\ga_* Y_P)= \Theta'(Y_{P'})=Y$.
\item[b)] $r_g^*\ga^*\Theta' = \ga^*r_g^*\Theta ' = \ga^*(\Ad g\inv)\Theta'=
(\Ad g\inv) \ga^*\Theta'$.
\end{enumerate}
Hence, $\ga^*\Theta'$ is a connection for $\pi: P\to M$.

\smallskip

\noindent
$(\Rightarrow)$ Suppose $\ga^*\Theta'$ is a connection for $\pi:P\to M$, 
and $Y\in \g$.  Then
\begin{enumerate}
\item[i)] $\Theta'(Y_{P'})= \Theta'(\ga_*Y_P)= (\ga^*\Theta')(Y_P)=Y$.
\item[ii)] $\ga^*r_g^*\Theta'= r_g^*\ga^* \Theta' = (\Ad g\inv) \ga^*\Theta'
=\ga^*((\Ad g\inv)\Theta ')$.
\end{enumerate}
Since $\ga$ is a submersion, $\ga^*$ is injective; hence,
\[
r_g^* \Theta' = (\Ad g\inv) \Theta'.
\]
So $\Theta'$ is a connection for $\pi': P' \to M'$. 
\medskip

Returning to the situation preceding the lemma, 
we will identify $\Xi'$ with $\Xi$.  This is possible because of the 
correspondence between forms on $P^{r}$ and basic forms on $\stl\times 
P$.

\separate
\begin{proposition} \label{31p:eqcurvature}
Let $K$ be the curvature of the connection $\Theta$ on the 
$\sone$-equivariant principal $G$-bundle $P\to M$.  Then the curvature $K_{\Xi}$ 
of the induced connection $\Xi$ on $P^{r} \to M^{r}$ is the 
following $S^1$-basic form on $\stl \times P$:
\[
K_{\Xi} = 1 \otimes K + d\theta \otimes L - \theta\otimes \ix K.
\]
\end{proposition}

\begin{proof}
\begin{align*}
d\Xi &= 1\otimes d\Theta + (d\theta)\otimes L - \theta \otimes dL, \\
[\Xi, \Xi] &= 
[1\otimes \Theta + \theta\otimes L, 1\otimes \Theta + \theta\otimes L]\\
&=1\otimes [\Theta,\Theta]+\theta\otimes [L,\Theta] -\theta\otimes 
[\Theta, L] \\
&= 1\otimes [\Theta,\Theta] -\theta\otimes 
2 [\Theta, L].
\end{align*}
To compute a bracket such as $[1\otimes \Theta, \theta\otimes L]$, first 
choose a basis $e_1, \ldots, e_n$ for $\g$ and write
\[
\Theta=\sum \Theta^i e_i, \quad L= \sum L^j e_j,
\]
with $\Theta^i, L^j$ being ordinary forms on $P$.  Then
\begin{align*}
[1\otimes \Theta,\, \theta\otimes L] &= \sum [ 1\otimes \Theta^i e_i,\, 
\theta\otimes L^j e_j] \\
&= \sum \Theta^i \theta L^j [e_i, e_j] \\
&=-\sum\theta \Theta^i L^j [e_i, e_j] \\
&=-\theta[\Theta^i e_i, L^j e_j] \\
&= -\theta[\Theta, L] \\
&=-\theta\otimes[\Theta, L].
\end{align*}
Therefore,
\begin{align*}
K_{\Xi} &= d\Xi + \frac{1}{2} [\Xi, \Xi]\\
&=1\otimes (d\Theta +\half [\Theta, \Theta])+(d\theta) \otimes L
-\theta\otimes(dL+[\Theta, L])\\
&= 1\otimes K +d\theta\otimes L - \theta \otimes DL \\
&= 1\otimes K +d\theta\otimes L - \theta \otimes \ix K \quad \mbox{(by 
Lemma \ref{p:ixk}).} 
\end{align*}

Since $\Xi$ is basic with respect to the $S^1$-action, it follows easily 
that $\iota_X K_{\Xi}=0$ and $\Lie_X K_{\Xi}=0$.  Hence, $K_{\Xi}$ is also basic
with respect to the $S^1$-action.
\end{proof}

\subsection{Equivariant Curvature and Equivariant Characteristic Classes}

In Proposition \ref{31p:eqcurvature}, the forms $\Xi$ and $\theta$ depend 
on $r$, but $K$ and $L$ do not.  To indicate this dependence, we write 
$K_{r}$ for $K_{\Xi}$ and $\theta_{r}$ for $\theta$.  
Let $u_r= d\theta_r$.
Then, dropping 
the tensor product sign, the curvature of the induced connection on 
$P^{r}\to M^{r}$ is
\[
K_{r} = K+u_{r}L - \theta_{r}(\ix (K+u_{r}L)) \in 
\Omega^*(P^{r}).
\]
Under the sequence of inclusions
\[
\begin{matrix}
P^1 &\to &P^2 &\to &\cdots & P^{r}& \to & \cdots & \to & P_{\sone} \\
\downarrow &&\downarrow &&&\downarrow &&&&\downarrow \\
M^1 &\to &M^2 &\to &\cdots & M^{r}& \to & \cdots & \to & M_{\sone}
\end{matrix},
\]
an element in the Weil model of $P_{S^1}$ that restricts to all 
$K_{r}$ is
\[
K_{\infty}:= K+uL-\theta (\ix (K+uL)),
\]
where $u$ and $\theta$ now have their usual meaning in the Weil model.
Moreover, by ((\cite{hatcher}), Prop.\ 3.52), in any dimension $k$ 
the real cohomology of the homotopy quotient 
$M_{S^1}$ is the inverse limit of $H^k(M^{r})$ as $r \to \infty$.

\separate
\begin{proposition}
If $f$ is an $\Ad (G)$-invariant polynomial on $\g$, then $f(K_{\infty}) 
\in W(\fraks) \otimes \Omega^*(P)$ is basic with respect to both the 
$\sone$-action and the $G$-action on $P$.
\end{proposition}

\begin{proof}
To check that $K_{\infty}$ is basic with respect to the $\sone$-action,
we compute for $X \in \Liealg (\sone)$:
\begin{align*}
\iota_X K_{\infty} &= \iota_X K + u \iota_X L - \iota_X(K+uL) 
=\iota_X K - \iota_X K =0 ,\\
\Lie_X K_{\infty} &= \Lie_X K + u \Lie_X L - \theta \Lie_X \iota_X (K+uL) \\
&\phantom{xxxxxxxxxxxxxxxxxxxx} \text{(in the circle case, $\Lie_X u = \Lie_X \theta = 0$)}\\
&= 0  \,\,\qquad \qquad \text{(since $\Lie_X K=0, \Lie_X L=0$, and 
$\Lie_X\iota_X=\iota_X \Lie_X$).}
\end{align*}

To check that $K_{\infty}$ is basic with respect to the $G$-action, let 
$Z\in \g$.  Since $K$ and $L$ are both horizontal with respect to $G$, 
$\iota_Z K_{\infty}=0$.

If $g\in G$, then by the $G$-equivariance of $K$ and $L$,
\[
r_g^* K_{\infty}= (\Ad g\inv) K_{\infty}.
\]
Hence,
\[
r_g^* f(K_{\infty})= f(r_g^* K_{\infty})= f((\Ad g\inv)K_{\infty})= 
f(K_{\infty}).
\]
Since $f(K_{\infty})$ is $G$-invariant, $L_Z K_{\infty}=0$ for any $Z \in 
\g$. 
\end{proof}

It follows from this proposition that $f(K_{\infty})$ is an element of 
$\Basic_{\sone} (W(\fraks)\otimes \Omega^*(P)) $ that descends to 
$\Basic_{\sone} (W(\fraks)\otimes \Omega^*(M))$.  As it is the unique 
element in the Weil model that restricts to $f(K_{r})$ for all $r$, 
its cohomology class represents the equivariant characteristic class 
$c_f(P)$.

Recall that the recipe for going from the Cartan model to the Weil model 
of the equivariant cohomology of a circle action is
\[
a \mapsto a-\theta (\ix a).
\]
This shows that the equivariant characteristic class of $P\to M$
 corresponding to an 
$\Ad (G)$-invariant polynomial $f$ on $\g$ is represented by $f(K+uL)$ in 
the Cartan model of $H_{\sone}^*(M)$.  In particular, 
$K+uL$ is the correct notion 
of \emph{equivariant curvature} in the Cartan model.

\section{The General Case}

Returning to the general case of a compact connected Lie group $S$, we 
approximate the universal $S$-bundle $ES \to BS$ by principal $S$-bundles 
of finite-dimensional manifolds
\[
\begin{matrix}
ES(1) &\hookrightarrow \cdots \hookrightarrow &ES(r) 
&\hookrightarrow & ES(r+1) & 
\hookrightarrow \cdots  \hookrightarrow & ES \\
\downarrow &&\downarrow &&\downarrow &&\downarrow \\
BS(1) &\hookrightarrow \cdots \hookrightarrow &BS(r) 
&\hookrightarrow& BS(r+1)& 
\hookrightarrow  \cdots  \hookrightarrow & BS
\end{matrix}.
\]
On each principal $S$-bundle $ES(r) \to BS(r)$, choose a connection 
$\theta(r)$.  Now fix an $r$ and for simplicity write $\theta$ instead of 
$\theta(r)$.

Let $\tilde{X}_1, \dots, \tilde{X}_{\ell}$ denote the fundamental vector fields 
on $ES(r)$ corresponding to the basis $\seq{X}{\ell}$ for the Lie algebra 
$\fraks$.  If $\theta=\sum\theta_i X_i$, then
\begin{enumerate}
\item[i)] $\theta (\tilde{X_j})= \sum \theta_i(\tilde{X}_j) X_i = X_j 
\quad \Rightarrow \quad \theta_i(\tilde{X}_j)= \delta_{ij}$.
\item[ii)] for $s\in S$, $r_s^* \theta = (\Ad s\inv) \theta$.
\end{enumerate}

As before, we adopt the shorthand  $\iota_i=\iota_{X_i}$, 
$\call_i=\call_{X_i}$, and $L_i= L_{X_i}$.  

\separate
\begin{proposition}
Let $c_{ij}^k$ be the structure constants of the Lie algebra $\fraks$ with 
respect to the basis $\seq{X}{\ell}$: $[X_i, X_j]= \sum c_{ij}^k X_k$.  
If $\theta=\sum \theta_i X_i$ is a connection on a principal $S$-bundle, 
then
\[
\call_i \theta_j = - \sum c_{ik}^j \theta_k.
\]
\end{proposition}

\begin{proof}
Since $\theta_j (\tilde{X}_k)= \delta_{jk}$,
\begin{align*}
0&= \call_i (  \theta_j (\tilde{X}_k)) \\
&= (\call_i \theta_j)(\tilde{X}_k)+ \theta_j (\call_i \tilde{X}_k).
\end{align*}
Hence,
\[
(\call_i \theta_j) \tilde{X}_k= - \theta_j ([\tilde{X}_i, \tilde{X}_k])
= -\theta_j (\sum c_{ik}^{\ga} \tilde{X}_{\ga})=-c_{ik}^j.
\]
If $\call_i \theta_j = \sum b_{ik}^j \theta_k$, then $b_{ik}^j=
(\call_i \theta_j) (\tilde{X}_k)=-c_{ik}^j$.
\end{proof}

Let $\Theta$ be an $S$-invariant connection on the $S$-equivariant 
principal $G$-bundle $\pi:P\to M$.  Define on $ES(r)\times P$ the 
$\g$-valued 1-form
\[
\Xi = 1 \otimes \Theta + \sum \theta_i \otimes L_i.
\] 

\separate
\begin{theorem}
The $\g$-valued 1-form $\Xi$ is a connection for the principal $G$-bundle
\[
ES(r) \times P \to ES(r) \times M.
\]
Moreover, $\Xi$ is basic with respect to the $S$-action on $ES(r)\times P$.
\end{theorem}

\begin{proof}
First we check that $\Xi$ is a connection on $ES(r)\times P$.
\begin{enumerate}
\item[i)] Let $Z \in \g$ and let $Z_P^*$ be its fundamental vector field 
on $P$.  Then the fundamental vector field of $Z$ on $ES(r)\times P$ is 
$(0,Z_P^*)$, because $G$ acts trivially on the first factor.  Hence,
\[
\Xi (0,Z_P^*) = 1\cdot \Theta(Z_P^*) +0 = Z.
\]
\item[ii)] Right-equivariance:
\begin{align*}
r_g^* \Xi &= 1\otimes r_g^*\Theta + \sum \theta_i \otimes r_g^* L_i = (\Ad 
g\inv) (1 \otimes \Theta + \sum \theta_i \otimes L_i) \\
&=(\Ad g\inv)\Xi.
\end{align*}
\end{enumerate}

Next we check that $\Xi$ is basic with respect to the $S$-action.
\begin{enumerate}
\item[iii)] Horizontality:
\begin{align*}
\iota_{X_j} \Xi &= 1\otimes \iota_{X_j}\Theta + \sum 
(\iota_{X_j}\theta_i)\otimes L_i + \sum \theta_i \otimes \iota_{X_j} L_i \\
&= -1 \otimes L_j + 1\otimes L_j + 0 = 0.
\end{align*}

\item[iv)] $S$-invariance:
\begin{align*}
\call_{X_j} \Xi &= 1\otimes \call_{X_j}\theta + \sum 
(\call_{X_j}\theta_i)\otimes L_i + \sum \theta_i \otimes \call_{X_j} L_i \\
&= 0+ \sum -c_{jk}^i \theta_k \otimes L_i + \sum \theta_i \otimes L_{[X_j, 
X_i]}  \qquad \text{(by Prop.~\ref{p:lielx})}\\
&= -\sum c_{jk}^i \theta_k \otimes L_i + \sum c_{ji}^k \theta_i \otimes 
L_k = 0.  \phantom{xxxxxxxxxxxxxx}
\end{align*}
\end{enumerate}
\end{proof}

Let $P^r = ES(r)\times_S P$ and $M^r=ES(r) \times_S M$.  The principal 
$G$-bundle $P^r\to M^r$ is our finite-dimensional approximation to the 
principal $G$-bundle of homotopy quotients $P_S \to M_S$.
By Lemma \ref{p:pullback}, $\Xi$ is the pullback of a connection $\Xi'$ 
on $P^r\to M^r$.  We will identify $\Xi$ with $\Xi'$.

\separate
\begin{proposition}
Let $K$ be the curvature of the $S$-invariant connection $\Theta$ on the 
$S$-equivariant principal $G$-bundle $P\to M$.  Then the curvature 
$K_{\Xi}$ of the induced connection $\Xi$ on $P^r \to M^r$ is given by the 
following $S$-basic form on $ES(r) \times P$:
\begin{equation} \label{e:kxi1}
K_{\Xi} = 1 \otimes K + \sum (d\theta_i) \otimes L_i - \sum \theta_i 
\otimes \iota_i K + \sum_{i<j} \theta_i\theta_j \otimes [L_i, L_j].
\end{equation}
\end{proposition}

\begin{proof}
\begin{align*}
d\Xi &= 1\otimes d\Theta + \sum (d\theta_i) \otimes L_i - \sum \theta_i 
\otimes dL_i.\\
[\Xi, \Xi] &= [1\otimes \Theta+\sum \theta_i \otimes L_i ,\ 1\otimes 
\Theta + \sum \theta_i \otimes L_i] \\
&= 1\otimes [\Theta, \Theta] + \sum \theta_i \otimes [L_i, \Theta] -
\sum \theta_i \otimes [\Theta, L_i]+ \sum_{i,j} \theta_i\theta_j \otimes
[L_i, L_j] \\
&= 1\otimes [\Theta, \Theta] -
\sum \theta_i \otimes 2 [\Theta, L_i]+ \sum_{i,j} \theta_i\theta_j \otimes
[L_i, L_j]
\end{align*}
In the last term, if we sum over only $i,j$ such that $i < j$, instead of 
all $i,j=1, \dots, \ell$, it becomes
\[
\sum_{i,j=1}^{\ell} \theta_i\theta_j \otimes [L_i, L_j]
=2 \sum_{i < j}\theta_i\theta_j \otimes [L_i, L_j].
\]
Hence,
\begin{align*}
K_{\Xi} &= d\Xi + \frac{1}{2}[\Xi, \Xi] \\
&=1\otimes K + \sum (d\theta_i) \otimes L_i - \sum \theta_i \otimes 
(dL_i + [\Theta, L_i]) + \sum_{i<j} \theta_i\theta_j \otimes [L_i, L_j].
\end{align*}
The final formula now follows from Prop.~\ref{p:ixk}(i).
\end{proof}

For simplicity, we drop the tensor product symbol in Eq.\ \eqref{e:kxi1}.  
Let 
\[
u_k= d\theta_k + \sum_{i<j} c_{ij}^k \theta_i\theta_j
\]
 be the curvature of 
the connection $\theta$ on $ES(r) \to BS(r)$.  We can rewrite \eqref{e:kxi1} as
\begin{align} \label{e:kxi2}
K_{\Xi} &= K+ \sum_k (d\theta_k + \sum_{i<j} c_{ij}^k \theta_i\theta_j) L_k
-\sum \theta_i \iota_i K \notag \\
&\phantom{xxxxxxxxxxxxxxxx}
 + \sum_{i<j} \theta_i \theta_j [L_i, L_j] - \sum_k \sum_{i<j} 
\theta_i\theta_j c_{ij}^k L_k \notag \\
&= K+ \sum u_kL_k - \sum \theta_i \iota_i K+ 
\sum_{i<j}\theta_i\theta_j[L_i, L_j] -\sum_i \sum_{i<j} \theta_i\theta_j 
c_{ij}^k L_k.
\end{align}

The last two terms of \eqref{e:kxi2} can be simplified somewhat, for
\[
\sum_k c_{ij}^k L_k = \sum c_{ij}^k L_{X_k} = L_{\sum c_{ij}^k X_k}
= L_{[X_i, X_j]};
\]
hence,
\begin{align*}
[L_i, L_j] - \sum_k c_{ij}^k L_k &= [L_{X_i}, L_{X_j}] - L_{[X_i,X_j]} \\
&= \iota_{X_j}\iota_{X_i} K \qquad\text{(by Prop.~\ref{p:ixk}(ii))}.
\end{align*}
Becaue $\iota_i$ anticommutes with $\iota_j$ and $\theta_j$, 
$\theta_i\theta_j\iota_j \iota_i = \theta_i \iota_i \theta_j \iota_j$
and
\[
K_{\Xi} = K+ \sum u_k L_k - \sum \theta_i \iota_i K+ \sum_{i 
<j}\theta_i\iota_i \theta_j\iota_j K.
\]
Since $L_k$ is a 0-form, $\iota_i L_k=0$ for all $i$.  Similarly, since $K$ is a 
2-form, contracting $K$ three or more times with vector fields yields 0, 
for example, $\iota_{i_1}\iota_{i_2}\iota_{i_3} K=0$.  Hence,
\begin{align} \label{e:kxi3}
K_{\Xi} &= (K+\sum u_kL_k) - \sum \theta_i \iota_i (K+\sum u_k L_k)\notag \\
&\phantom{xxxxxxxxxx}
+\sum \theta_i\iota_i\theta_j \iota_j (K+\sum u_kL_k) - \dots \notag \\
&=\prod_{i=1}^{\ell} (1-\theta_i\iota_i)(K+\sum u_kL_k).
\end{align}

As in the circle case, in this formula $K_{\Xi}$, $\theta_i$ and $u_k$ are 
all differential forms on $P^r$ and should be more properly written as
$K_{\Xi}(r)$, $\theta_i(r)$ and $u_k(r)$ to indicate their dependence on 
$r$.  If $f$ is an $\Ad G$-invariant polynomial on $\g$, then 
$f(K_{\Xi}(r))$ represents the characteristic class $c_f(P^r)$.  Clearly, 
the element in the Weil model of $P_S$ that restricts to all 
$f(K_{\Xi}(r))$ is $f(K_{\infty})$, where $K_{\infty}$ is given by the same 
formula as \eqref{e:kxi3},
\[ 
K_{\infty}=\prod_{i=1}^{\ell} (1-\theta_i\iota_i)(K+\sum u_kL_k),
\]
but $\theta_i$ and $u_k$ now have their usual meaning as elements of the 
Weil algebra of $S$.

The rest of the argument proceeds as in the circle case.  The upshot is 
that $f(K_{\infty})$ is the element in the Weil model of $H_S^*(M)$ which 
restricts to the characteristic classes $c_f(P^r)$ for all $r$, and therefore 
represents the equivariant characteristic class $c_f(P)$.  Under the 
Mathai-Quillen isomorphism, the corresponding element in the Cartan model 
is $f(K+\sum u_k L_k)$.

An element of the Cartan model $(S^*(\fraks^*)\otimes \Omega^*(M))^S$ may 
be viewed as a polynomial function $h: \fraks \to \Omega^*(M)$ which is 
$S$-equivariant in the following sense:
\[
h((\Ad s)X) = (s\inv)^*h(X).
\]
From this point of view, $f(K+\sum u_kL_k)$ is the function$: \fraks \to 
\Omega^*(M)$ given by
\begin{equation} \label{e:function}
X \mapsto f(K+L_X) \quad \in \Omega^*(M),
\end{equation}
since
\[
L_X = L_{\sum u_k(X)X_k} = \sum u_k(X)L_{X_k} = \sum u_k L_k (X).
\]
By Prop.~\ref{p:translatelx},
\[
L_{(\Ad s)X} = (s\inv)^*L_X,
\]
which shows that $K+L_X$ is $S$-invariant and therefore $f(K+L_X)$ is an 
$S$-invariant element of $S^*(\fraks^*)\otimes \Omega^*(M)$.  Finally, if 
we denote the function in \eqref{e:function} by $K+L$, then the element of the 
Cartan model of $M_S$ corresponding to the equivariant 
characteristic class $c_f(P_S)$ 
is $f(K+L)$.

\end{document}